\documentclass[reqno,a4paper,12pt]{amsart} 

\usepackage{amsmath}
\usepackage{amscd,amsfonts,amssymb}
\usepackage{mathrsfs,dsfont}
\usepackage{pgfplots}
\pgfplotsset{compat=1.15}
\usepackage{tikz}
\usetikzlibrary{patterns}
\usepackage{float}
\usepackage{appendix}
\usepackage{caption}
\usepackage{subcaption}
\usepackage{enumerate}
\usepackage{hyperref}

\numberwithin{equation}{section}
\numberwithin{figure}{section}

\def\R{\mathbb{R}}

\def\1{\mathds{1}}

\renewcommand\leq{\leqslant}
\renewcommand\geq{\geqslant}

\renewcommand\hat{\widehat}
\renewcommand\Re{\operatorname{Re}}

\newcommand{\supp}{\operatorname{supp}}

\newcommand{\diam}{\operatorname{diam}}

\theoremstyle{plain}
\newtheorem{thm}{Theorem}[section]

\newtheorem{corollary}[thm]{Corollary}

\newtheorem*{claim*}{Claim}
\newtheorem*{thm*}{Theorem}

\theoremstyle{definition}

\newtheorem*{definition*}{Definition}
\newtheorem*{remarks*}{Remarks}
\newtheorem*{remark*}{Remark}

\newenvironment{enumerate-math}
{\begin{enumerate}
\addtolength{\itemsep}{5pt}
}
{\end{enumerate}}

\newenvironment{enumerate-text}
{\begin{enumerate}
\addtolength{\itemsep}{5pt}
}
{\end{enumerate}}

\begin{document}

\title{Fourier frames on smooth surfaces with nonvanishing Gaussian curvature}

\author{Xinyu Chen}
\email{12432007@mail.sustech.edu.cn}
\author{Bochen Liu}
\email{Bochen.Liu1989@gmail.com}
\address{Department of Mathematics \& International Center for Mathematics, Southern University of Science and Technology, Shenzhen 518055, China}

\thanks{This work was partially supported by the National Key R\&D Program of China 2024YFA1015400, and the National Natural Science Foundation of China grant 12131011.}
\date{}


\begin{abstract}
It is known that a small spherical cap (rigorously its surface measure) admits Fourier frames, while the whole sphere does not. In this paper, we prove more general results. Consequences indclude that a small spherical cap in $\mathbb{R}^d$ near the north pole cannot have a frame spectrum near the $x_d$-axis, and $S$ does not admit any Fourier frame if its interior contains a closed hemisphere. We also resolve the endpoint case, that is, a hemisphere does not admit any Fourier frame. This answers a question of Kolountzakis and Lai. Our results also hold on more general smooth surfaces with nonvanishing Gaussian curvature. In particular, any compact $(d-1)$-dimensional smooth submanifold immersed in $\mathbb{R}^d$ with nonvanishing Gaussian curvature does not admit any Fourier frame. This generalizes a previous result of Iosevich, Lai, Wyman and the second author on the boundary of convex bodies, as well as improves a recent result of Kolountzakis and Lai from tight frame to frame.

\end{abstract}
\maketitle

\section{Introduction}
Suppose $\mu$ is a finite Borel measure on $\R^d$. A discrete set $\Lambda\subset\R^d$ is called a frame
spectrum of $\mu$ if there exist $0<A\leq B<\infty$ such that
$$A\|f\|_{L^2(\mu)}^2\leq \sum_{\lambda\in\Lambda}|\hat{f\,d\mu}(\lambda)|^2\leq B\|f\|_{L^2(\mu)}^2,\ \forall\,f\in L^2(\mu).$$
In this case we say $\mu$ admits a Fourier frame. It is called a tight frame if $A=B$.

Fourier frame was introduced by Duffin and Schaeffer \cite{DS52} in 1952, and has become a powerful tool in applied harmonic analysis. It is known \cite{NOU16} that any Euclidean subset of finite Lebesgue measure admits Fourier frames (with $d\mu=\chi_Edx$). Because of this strong result, the discussion on the existence of Fourier frames mainly focuses on singular measures. See \cite{Lai11}\cite{HLL13}\cite{DL14}\cite{FL18}\cite{LL25+} for examples admitting no Fourier frames.

We say a subset of a surface admits Fourier frames if the restricted surface measure on it does.

In \cite{Lev18}, Lev pointed out that a small spherical cap admits Fourier frames, and then asked about the whole sphere. Lev's question was answered negatively by Iosevich, Lai, Wyman and the second author in \cite{ILLW19}. The proof relies on the well-known asymptotic formula of the Fourier transform of the surface measure on the sphere:
$$\widehat{\sigma}(\xi) = C\left(\frac{\xi}{|\xi|}\right)|\xi|^{-\frac{d-1}{2}} \cos\left( 2\pi \left(|\xi|-\frac{d-1}{8}\right)\right) + O(|\xi|^{-\frac{d-1}{2}-1}).$$
Consequently,
\begin{equation}\label{Fourier-upper-bound}|\hat{\sigma}(\xi)|\lesssim |\xi|^{-\frac{d-1}{2}}	
\end{equation}
and
\begin{equation}\label{Fourier-lower-bound}
\int_{B_1(\xi)}|\hat{\sigma}(\eta)|^2d\eta\gtrsim |\xi|^{-(d-1)}\int_{|\xi|-1}^{|\xi|+1} \left|\cos(2\pi(t-\frac{d-1}{8})\right|^2dt\gtrsim |\xi|^{-(d-1)}.
\end{equation}
More generally, the same conclusion and estimates \eqref{Fourier-upper-bound}\eqref{Fourier-lower-bound} hold on the smooth boundary of any convex body with nonvanishing Gaussian curvature, by the famous Herz formula \cite{Herz62}. As a comparison, the boundary of any polytope in $\R^d$ admits Fourier frames \cite{ILLW19}, while in the plane these frames can never be tight according to a recent work of Kolountzakis and Lai \cite{KL25+}. They also considered self-intersecting curves, showing that a smooth closed planar curve admits no tight Fourier frames if it has positive curvature everywhere and there are only finitely many self-intersections all of which are transverse. See Theorem 2 in \cite{KL25+}.

In this paper, we prove more general results on more general smooth surfaces with nonvanishing Gaussian curvature. Let $U_i\subset\R^{d-1}$, $d\geq 2$, $i=1,\dots,M$, be open sets, and $\chi_i: U_i\rightarrow\R^d$, be smooth embeddings such that each surface $\chi_i(U_i)$ has nonvanishing Gaussian curvature, with $\sigma_i$ as the surface measure. Then let
\begin{equation}\label{def-S}S:=\bigcup_i\chi_i(U_i)\end{equation}
be a smooth $(d-1)$-dimensional submanifold immersed in $\R^d$ and define its surface measure by
\begin{equation}\label{def-surface-measure-S}d\sigma_S(p):=\sum_i m(p)^{-1}d\sigma_i(p),\end{equation}
where $m(p)$ is the multiplicity function 
$$m(p):=\#\{i:p\in \chi_i(U_i)\}.$$
Notice that $\sigma_S$ only depends on $S$, not on the choice of $\chi_i(U_i)$. Given a point $p\in\chi_i(U_i)$, let $\mathcal{C}_p$ denote the $1$-dimensional subspace along its normal, and 
\begin{equation}\label{def-C-S}\mathcal{C}_S:=\bigcup_i\bigcup_{p\in\chi_i(U_i)}\mathcal{C}_{p},\end{equation}
that is also independent in the choice of $\chi_i(U_i)$.

We say $S'$ is a $(d-1)$-dimensional submanifold compactly contained in $S$, if there exist open sets $V_i\subset U_i$ (not necessarily nonempty) with closure $\overline{V_i}\subset U_i$ such that
\begin{equation}\label{def-S-prime}S'=\bigcup_i\chi_i(V_i).\end{equation}
Similarly, one can define 
\begin{equation}\label{def-C-S-prime}\mathcal{C}_{S'}:=\bigcup_i\bigcup_{p\in\chi_i(V_i)}\mathcal{C}_{p}:=\bigcup_i\mathcal{C}_i.\end{equation}

\begin{thm}\label{thm-cap}
	Let $S, \sigma_S$ be as above. Suppose $\sigma_S$ is finite and admits a frame spectrum $\Lambda$. Then
	$$\Lambda\cap (\mathcal{C}_{S'})^c\neq\emptyset$$
    for every $(d-1)$-dimensional submanifold $S'$ compactly contained in $S$. In particular, if there exists such an $S'$ with $\mathcal{C}_{S'}=\R^d$, then $\sigma_S$ admits no Fourier frames.
\end{thm}
We need $S'\subset S$ to state Theorem \ref{thm-cap} due to technical reasons. If $S$ is a compact subset in $\R^d$, then a simple topological argument shows that $S'=S$ is allowed. More precisely, for each point $p\in S$, there exists some $i_p$ and an open set $W_p$ satisfying $$\chi_{i_p}^{-1}(p)\in W_p\subset\overline{W_p}\subset U_{i_p}.$$
Then the compactness guarantees that we only need finitely many such $p_j$ for $$S=\bigcup_j\chi_{i_{p_j}}(W_{p_j}),$$and $S'$ can be chosen to be equal to $S$ by taking
$$V_i:=\bigcup_{j:W_{p_j}\subset U_i}W_{p_j}.$$
Moreover, in this case, for each $e\in S^{d-1}$, by translating hyperplanes perpendicular to $e$ for tangency at $S$, one can always find some $i$ and some $p\in\chi_{i}(U_{i})$ with $e$ as a normal. Hence Theorem \ref{thm-cap} implies the following. 
\begin{corollary}\label{cor}
    Let $S$ be a compact $(d-1)$-dimensional smooth submanifold immersed in $\R^d$ with nonvanishing Gaussian curvature defined in \eqref{def-S}, then its surface measure $\sigma_S$ defined in \eqref{def-surface-measure-S} does not admit any Fourier frame.
\end{corollary} 
Corollary \ref{cor} generalizes the result of Iosevich, Lai, Wyman and the second author in \cite{ILLW19} on the smooth boundary of convex bodies with nonvanishing Gaussian curvature. One may wonder what a surface in Corollary \ref{cor} may look like if it is not the smooth boundary of a convex body. Kolountzakis and Lai drew an example of self-intersecting curve in the plane (see Figure 3 in \cite{KL25+}). There also exist explicit parameterizations such as 
$$(a(\theta),b(\theta)):=(\cos\theta+2\cos 2\theta, \sin\theta+\sin 2\theta), \ 0\leq \theta\leq 2\pi.$$
See Figure \ref{Figure-1} below. The previous result in \cite{ILLW19} does not apply to such self-intersecting curves due to the lack of convexity, while by Corollary \ref{cor} they do not admit Fourier frames. This also improves a result of Kolountzakis and Lai (Theorem 2 in \cite{KL25+}) from tight frame to frame.

\begin{figure}[h]
\centering
\begin{tikzpicture}

\begin{axis}[
    width=9cm,
    height=6cm,
    xlabel={$x$},
    ylabel={$y$},
    axis lines=middle,
    xmin=-3, xmax=3,
    ymin=-2, ymax=2,
    grid=both,
]

\addplot [
    domain=0:2*pi,
    samples=500,
    smooth,
    thick,
    blue
] 
(
    {cos(deg(x)) + 2*cos(2*deg(x))}, 
    {sin(deg(x)) + sin(2*deg(x))}
);

\end{axis}
\end{tikzpicture}
\caption{$(x,y)=(\cos\theta + 2\cos 2\theta,  \sin\theta + \sin 2\theta)$}\label{Figure-1}
\end{figure}
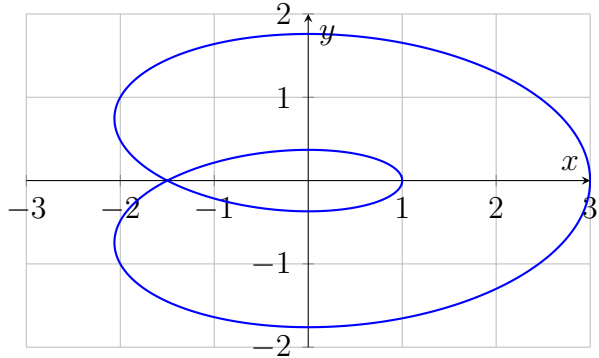
\begin{figure}[h]
\centering
\begin{tikzpicture}

\begin{axis}[
    width=9cm,
    height=6cm,
    xlabel={$x$},
    ylabel={$y$},
    axis lines=middle,
    xmin=-3, xmax=3,
    ymin=-2, ymax=2,
    grid=both,
]

\addplot [
    domain=0:2*pi,
    samples=500,
    smooth,
    thick,
    blue
] 
(
    {cos(deg(x)) + 2*cos(2*deg(x))}, 
    {sin(deg(x)) + sin(2*deg(x))}
);
\addplot [
    domain=0:2*pi,
    samples=500,
    smooth,
    thick,
    blue
] 
(
    {-cos(deg(x)) - 2*cos(2*deg(x))}, 
    {sin(deg(x)) + sin(2*deg(x))}
);

\end{axis}
\end{tikzpicture}
\caption{It is the slice of $(a(\theta)\cos\varphi, a(\theta)\sin\varphi, b(\theta))$ in the $xy$-plane. The curve that can be observed from outside is a reflection of the right half of Figure \ref{Figure-1}, that is not $C^3$ at the north pole and the south pole.}\label{Figure_2}
\end{figure}
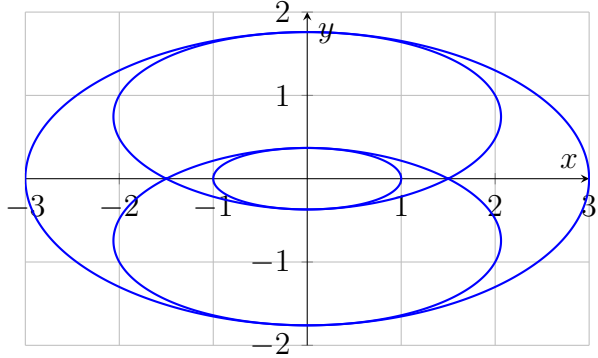

In $\R^3$, one can rotate this curve in the $xy$-plane about the $y$-axis to obtain
$$(a(\theta)\cos\varphi, a(\theta)\sin\varphi, b(\theta)), \ 0\leq \theta, \varphi\leq 2\pi,$$
on which, thanks to $b'=a$, the Gaussian curvature equals 
$$\frac{b''(\theta)a'(\theta)-a''(\theta)b'(\theta)}{(a'(\theta)^2+b'(\theta)^2)^2}=\frac{17+6\cos\theta+4\cos^3\theta}{(a'(\theta)^2+a(\theta)^2)^2}\neq 0.$$
Although the rotation makes it look like an ellipsoid from outside, the boundary is not $C^3$. See Figure \ref{Figure_2} below for an explanation. Therefore the formula of Herz does not apply (see Theorem 3 in \cite{Herz62}). Consequently one cannot apply the previous result in \cite{ILLW19} to the ``boundary" observed from outside to conclude the nonexistence of Fourier frames. On the other hand, our Corollary \ref{cor} applies.

Similar examples in general $\R^d$, $d\geq 3$, can be taken as the revolution surface 
$$(a(\theta)\omega, b(\theta)), \ \omega\in S^{d-2}, 0\leq \theta\leq 2\pi,$$
on which the Gaussian curvature equals 
$$\frac{b''(\theta)a'(\theta)-a''(\theta)b'(\theta)}{(a'(\theta)^2+b'(\theta)^2)^{\frac{d+1}{2}}}=\frac{17+6\cos\theta+4\cos^3\theta}{(a'(\theta)^2+a(\theta)^2)^{\frac{d+1}{2}}}\neq 0.$$

On the sphere, Theorem \ref{thm-cap} implies that there is no way for a frame spectrum for a spherical cap near the north pole to lie near the $x_d$-axis.  It also implies that, if the interior of a spherical cap contains a closed hemisphere, then it does not admit any Fourier frame. Recall that Lev pointed out that a spherical cap admits Fourier frames if it is compactly contained in a hemisphere. So the hemisphere is the threshold, but not covered by any previous result including our Theorem \ref{thm-cap}. This question is also asked by Kolountzakis and Lai recently (see Question 2 at the end of \cite{KL25+}).
\begin{thm}\label{thm-hemi}
	Suppose $S$ is the boundary of a centrally symmetric convex body in $\R^d$ centered at the origin, smooth, and has nonvanishing Gaussian curvature everywhere. 	Then the restricted surface measure $\sigma_+$ on $S_+:= S\cap \{x_d\geq 0\}$ does not admit any Fourier frame.
\end{thm}
As a remark, although the proof on the hemisphere is inspired by the whole sphere, the surface measure $\sigma_+$ satisfies different Fourier properties from $\sigma$ when $d\geq 4$: by the spherical coordinate and integration by parts,
\begin{equation}\label{lower-Fourier-bound-hemisphere}|\widehat{\sigma_+}(0,\dots,0,\xi_d)|=C_d\left|\int_0^{\pi/2}e^{-2\pi i \xi_d\cos\theta}(\sin\theta)^{d-2}d\theta\right|\gtrsim |\xi_d|^{-1}.\end{equation}
In particular, $\sigma_+$ is not a Salem measure when $d\geq 4$. On the other hand, $S^{d-1}_+$ and $S^{d-1}$ share the same Fourier properties on $\R^{d-1}\times\{0\}$, because 
$$\widehat{\sigma_+}(\xi_1,\dots,\xi_{d-1}, 0)=\frac{1}{2}\widehat{\sigma}(\xi_1,\dots,\xi_{d-1},0)$$
by the symmetry on the last coordinate. In this paper, we find simple arguments that do not require a deep understanding of $\widehat{\sigma_+}$, but it would be interesting to figure out the Fourier property of the hemisphere, which seems not covered in the existing literature. We refer to Question 416545 in MathOverFlow for a discussion.

\subsection*{Notation.}

Throughout this article, we write $A\lesssim B$ if there exists an absolute constant $C$ such that $A\leq CB$.

\subsection*{Acknowledgments.}
We would like to thank Chun-Kit Lai for comments on the manuscript and discussion on his recent work joint with Kolountzakis.

\section{An observation on previous work and the proof of Theorem \ref{thm-hemi}}
In \cite{ILLW19} (see Theorem 1.3, 1.4 there), it is proved that \eqref{Fourier-upper-bound} together with 
\begin{equation}\label{lower-frame-bound}A\|f\|_{L^2(\sigma)}^2\leq \sum_{\lambda\in\Lambda}|\hat{f\,d\sigma}(\lambda)|^2,\ \forall\,f\in L^2(\sigma)\end{equation}
imply
\begin{equation}\label{divergence}\sum_{\lambda\in\Lambda\backslash \{0\}}|\lambda|^{-(d-1)}=\infty,\end{equation}
while \eqref{Fourier-lower-bound} together with 
\begin{equation}\label{upper-frame-bound}\sum_{\lambda\in\Lambda}|\hat{f\,d\sigma}(\lambda)|^2\leq B\|f\|_{L^2(\sigma)}^2,\ \forall\,f\in L^2(\sigma)\end{equation}
imply
\begin{equation}\label{convergence}\sum_{\lambda\in\Lambda\backslash \{0\}}|\lambda|^{-(d-1)}<\infty.\end{equation}

An observation that helps us in this paper is that, the proof in \cite{ILLW19} does not take use of all $f\in L^2$ from \eqref{lower-frame-bound}\eqref{upper-frame-bound}. In fact, it only considers the family $\{e^{2\pi i x\cdot\xi}\}_{\xi\in\R^d}$. More precisely,  \eqref{Fourier-upper-bound} together with
$$1\lesssim\sum_{\lambda\in\Lambda}|\hat{\sigma}(\lambda-\xi)|^2,\ \forall\,\xi\in\R^d,$$
are sufficient to conclude the divergence \eqref{divergence}. Also \eqref{Fourier-lower-bound} together with
$$\sum_{\lambda\in\Lambda}|\hat{\sigma}(\lambda-\xi)|^2\lesssim 1,\ \forall\,\xi\in\R^d,$$
are sufficient to conclude the convergence \eqref{convergence}.

With this observation in mind, Theorem \ref{thm-hemi} quickly follows, even if $\sigma_+$ has Fourier properties different from $\sigma$. 

Take $d\mu:=\psi d\sigma_+$, where $\psi\in C_0^\infty(S_+)$ is fixed. Then the stationary phase estimate (see, e.g. Theorem 1.2.1 in \cite{Sogge17}) implies that the upper Fourier bound \eqref{Fourier-upper-bound} holds on $\mu$. If $\sigma_+$ admits a frame spectrum $\Lambda$, with the family $f(x)=\psi(x)e^{2\pi i x\cdot\xi}$, $\forall\,\xi$, the lower frame bound \eqref{lower-frame-bound} implies
$$1\lesssim\sum_{\lambda\in\Lambda}|\hat{\mu}(\lambda-\xi)|^2,\ \forall\,\xi\in\R^d.$$
Hence the divergence \eqref{divergence} follows from our observation above.

On the other hand, by the central symmetry,
\begin{equation}\label{symmetry}2\Re\hat{\sigma_+}(\xi)=2\int_{S_+}\cos(2\pi i x\cdot\xi)d\sigma(x)=\int_{S}\cos(2\pi i x\cdot\xi)d\sigma(x)= \hat{\sigma}(\xi),\end{equation}
so $2|\hat{\sigma_+}|\geq |\hat{\sigma}|$ and the lower Fourier bound \eqref{Fourier-lower-bound} holds on $\sigma_+$. It is not the best estimate from below (recall \eqref{lower-Fourier-bound-hemisphere}), but enough for our use. As the upper frame bound \eqref{upper-frame-bound} with the family $\{e^{2\pi i x\cdot\xi}\}_{\xi\in\R^d}$ implies
$$\sum_{\lambda\in\Lambda}|\hat{\sigma_+}(\lambda-\xi)|^2\lesssim 1,\ \forall\,\xi\in\R^d,$$
the convergence \eqref{convergence} follows from our observation above, a contradiction.

\section{Proof of Theorem \ref{thm-cap}}
By our definition of $S$ in \eqref{def-S}, there exists $1\leq i_0\leq M$ and a nonzero $\psi\in C_0^\infty(\chi_{i_0}(U_{i_0}))$ on $S$. Also, by our definition of $\sigma_S$ in \eqref{def-surface-measure-S}, we have
$$\widehat{\psi d\sigma_S}=\widehat{\psi d\sigma_{i_0}}.$$
Then, by the stationary phase estimate (see, e.g. Theorem 1.2.1 in \cite{Sogge17}), the upper Fourier bound \eqref{Fourier-upper-bound} holds on $\psi d\sigma_S$. Therefore, if $\sigma_S$ admits a frame spectrum $\Lambda$, the divergence
$$\sum_{\lambda\in\Lambda\backslash\{0\}}|\lambda|^{-(d-1)}=\infty$$
follows in the same way as in the previous section.

It remains to show that, assuming $\Lambda\subset\mathcal{C}_{S'}$, then
$$\sum_{\lambda\in\Lambda\backslash\{0\}}|\lambda|^{-(d-1)}<\infty.$$
By our definition of $\mathcal{C}_{S'}$ in \eqref{def-C-S-prime}, it suffices to show
\begin{equation}\label{sum-C-i}\sum_{\lambda\in\mathcal{C}_i\cap  \Lambda\backslash\{0\}}|\lambda|^{-(d-1)}<\infty,\ \forall\,1\leq i\leq M.\end{equation}

Fix $i$ and recall our definition of $S'$ in \eqref{def-S-prime}. After translation and rotation, every $p\in \chi_i(\overline{V_i})$ has a $\delta_p$- neighborhood in $\chi_i(U_i)$ that becomes the graph $$(x_1,\dots,x_d, h_p(x_1,\dots, x_d))$$
centered at the origin with $h_p(0)=0$, $\nabla h_p(0)=0$. Moreover, by the Morse Lemma (see, e.g. Chapter 6 in \cite{Wol03}), this $\delta_p$-neighborhood can be taken such that 
$$h_p\circ G_p(y_1,\dots, y_{d-1})=\sum_{j=1}^k y_j-\sum_{j=k+1}^{d-1} y_j$$
under some smooth diffeomorphism $G_p$. Since $\overline{V_i}$ is compact, there exists a Lebesgue number $\delta>0$ associated with $\{B(p,\delta_p)\}_{p\in \chi_i(\overline{V_i})}$. Then there exist 
$$\psi_{ij}\in C_0^\infty(\chi_i(U_i)),\  j=1,\dots, N_i,$$ such that $\diam(\supp \psi_{ij})<\delta/4$ and $\sum_{j=1}^{N_i} \psi_{ij}=1$ on $\chi_i(\overline{V_i})$. In particular, for each $p\in \chi_i(\overline{V_i})$, there exists $j_0$ such that $\psi_{ij_0}(p)\geq \frac{1}{N_i}$. This means, to show \eqref{sum-C-i}, it suffices to show 
\begin{equation}\label{sum-C-ij}\sum_{\lambda\in\mathcal{C}_{ij}\cap  \Lambda\backslash\{0\}}|\lambda|^{-(d-1)}<\infty,\ \forall\,1\leq j\leq N_i,\end{equation}
where 
\begin{equation}\label{def-C-ij}\mathcal{C}_{ij}:=\mathcal{C}_i\cap \mathcal{C}_{\{\psi_{ij}\geq N_i^{-1}\}}.\end{equation}
Finally, by the upper frame bound \eqref{upper-frame-bound} with $f=\psi_{ij}$ and the fact $\widehat{\psi_{ij}d\sigma_{S}}=\widehat{\psi_{ij}d\sigma_i}$ (recall \eqref{def-surface-measure-S}), to prove \eqref{sum-C-ij} it suffices to show
    \begin{equation}\label{eq-lowerbound}
\vert\widehat{\psi_{ij}d\sigma_i}(\lambda)\vert^{2}\gtrsim \vert\lambda\vert^{-(d-1)},\ \forall\,\lambda\in\mathcal{C}_{ij}\backslash\{0\}.
    \end{equation}
    
    For each $\lambda\in\mathcal{C}_{ij}\backslash\{0\}$, there exists a unique $p_{\lambda}\in \supp\psi_{ij}$ with $\lambda$ as a normal. 
    Near this $p_{\lambda}$, the surface is locally the graph of $h_{p_{\lambda}}:=h_{\lambda}$ with $h_{\lambda}(0)=0$, $\nabla h_{\lambda}(0)=0$, after translation and rotation. More precisely, there exists a rotation $g_\lambda$ such that $$g_\lambda\lambda=(0,\dots,0,|\lambda|)$$
    and
    \begin{equation}\label{rotation-translation}\begin{aligned}
        \widehat{\psi_{ij} d\sigma_i}(\lambda)=&\int e^{-2\pi i (p_{\lambda}+g_\lambda^t(x, h_\lambda(x)))\cdot \lambda}\psi_{ij}(x,h_\lambda(x))\sqrt{1+|\nabla h_\lambda(x)|^2}\,dx\\=&e^{-2\pi i p_\lambda\cdot\lambda}\int e^{-2\pi i |\lambda|h_\lambda(x)}\psi_{ij}(x,h_\lambda(x))\sqrt{1+|\nabla h_\lambda(x)|^2}\,dx.
    \end{aligned}\end{equation}
    This is equal to, by the stationary phase estimate (see, e.g. Proposition 6.4 in \cite{Wol03}),
    \begin{equation}\label{error-term}e^{-2\pi ip_\lambda\cdot\lambda}e^{-\pi i\frac{s_\lambda}{4}} 2^{-\frac{d-1}{2}}|\det H_{h_\lambda}(0)|^{-\frac{1}{2}}\psi_{ij}(p_\lambda)\cdot|\lambda|^{-\frac{d-1}{2}}+O(|\lambda|^{-\frac{d-1}{2}-1}),\end{equation}
where $H_{h_\lambda}$ denotes the Hessian matrix of $h_\lambda$ and $s_{\lambda}$ denotes the signature of $H_{h_\lambda}(0)$. Here the implicit constant in $O(|\lambda|^{-\frac{d-1}{2}-1})$ may depend on derivatives of $h_\lambda$, but is uniform over $p_\lambda\in \chi_i(\overline{V_i})$ because locally $h_p$ can be chosen to be smooth in $p$. Also $|\det H_{h_\lambda}(0)|\gtrsim 1$ because $\det H_{h_\lambda}(0)$ is the Gaussian curvature at $p_\lambda\in \chi_i(\overline{V_i})$, and recall that $\psi_{ij}(p_{\lambda})\geq\frac{1}{N_i}$ for each $\lambda\in\mathcal{C}_{ij}$ by our construction of $\mathcal{C}_{ij}$ in \eqref{def-C-ij}. 

Overall, \eqref{rotation-translation} and \eqref{error-term} give
\[
    \vert\widehat{\psi_{ij}d\sigma_i}(\lambda)\vert^{2}\gtrsim\vert\lambda\vert^{-(d-1)},\ \forall\,\lambda\in\mathcal{C}_{ij}\backslash\{0\},
\]
that proves \eqref{eq-lowerbound}, thus completes the proof of Theorem \ref{thm-cap}.

As a final remark, in this section we do not follow the previous approach to pursue
\begin{equation}\label{bad-method}\int_{B_1(\xi)}|\hat{\sigma_S}(\eta)|^2d\eta\gtrsim |\xi|^{-(d-1)} \ \text{or}\ \int_{B_1(\xi)}|\hat{\psi d\sigma_S}(\eta)|^2d\eta\gtrsim |\xi|^{-(d-1)}, \quad \forall\,\xi\in\mathcal{C}_{S'}.\end{equation}
In fact we do not even need to take an integral over a neighborhood of $\xi$. Although \eqref{bad-method} is likely to be true, the argument would become more complicated than the above. By the stationary phase estimate (see, e.g. Proposition 6.4 in \cite{Wol03}), the value of the Fourier transform at $\eta$ is determined by points $p_\eta^k$ in the surface with $\pm\frac{\eta}{|\eta|}$ as a normal. As each normal vector may correspond to multiple points in $S$, the idea of considering \eqref{bad-method} ends up with an integral of form
$$\int_{B_1(\xi)}\left|\sum_k e^{-2\pi ip^k_{\eta}\cdot\eta}e^{-\pi i\frac{s^k_\eta}{4}} 2^{-\frac{d-1}{2}}|\det H_{h^k_\eta}(0)|^{-\frac{1}{2}}\psi_{ij}(p^k_{\eta})\right|^2d\eta.$$
Our ``local" argument in this paper successfully avoid dealing with such a complicated integral of exponential sums, that may shed light on other problems in the future.


\bibliographystyle{abbrv}
\bibliography{mybibtex.bib}

\begin{thebibliography}{10}

\bibitem{DS52}
R.~J. Duffin and A.~C. Schaeffer.
\newblock A class of nonharmonic {F}ourier series.
\newblock {\em Trans. Amer. Math. Soc.}, 72:341--366, 1952.

\bibitem{DL14}
D.~E. Dutkay and C.-K. Lai.
\newblock Uniformity of measures with {F}ourier frames.
\newblock {\em Adv. Math.}, 252:684--707, 2014.

\bibitem{FL18}
X.~Fu and C.-K. Lai.
\newblock Translational absolute continuity and {F}ourier frames on a sum of
  singular measures.
\newblock {\em J. Funct. Anal.}, 274(9):2477--2498, 2018.

\bibitem{HLL13}
X.-G. He, C.-K. Lai, and K.-S. Lau.
\newblock Exponential spectra in {$L^2(\mu)$}.
\newblock {\em Appl. Comput. Harmon. Anal.}, 34(3):327--338, 2013.

\bibitem{Herz62}
C.~S. Herz.
\newblock Fourier transforms related to convex sets.
\newblock {\em Ann. of Math. (2)}, 75:81--92, 1962.

\bibitem{ILLW19}
A.~Iosevich, C.-K. Lai, B.~Liu, and E.~Wyman.
\newblock Fourier frames for surface-carried measures.
\newblock {\em Int. Math. Res. Not. IMRN}, (3):1644--1665, 2022.

\bibitem{KL25+}
M.~N. Kolountzakis and C.-K. Lai.
\newblock Non-spectrality of some piecewise smooth curves and unions of line
  segments.
\newblock {\em arXiv:2507.00581}.

\bibitem{Lai11}
C.-K. Lai.
\newblock On {F}ourier frame of absolutely continuous measures.
\newblock {\em J. Funct. Anal.}, 261(10):2877--2889, 2011.

\bibitem{Lev18}
N.~Lev.
\newblock Fourier frames for singular measures and pure type phenomena.
\newblock {\em Proc. Amer. Math. Soc.}, 146(7):2883--2896, 2018.

\bibitem{LL25+}
L.~Li and B.~Liu.
\newblock Fourier frames on salem measures.
\newblock {\em arXiv preprint arXiv:2506.01280}.

\bibitem{NOU16}
S.~Nitzan, A.~Olevskii, and A.~Ulanovskii.
\newblock Exponential frames on unbounded sets.
\newblock {\em Proc. Amer. Math. Soc.}, 144(1):109--118, 2016.

\bibitem{Sogge17}
C.~D. Sogge.
\newblock {\em Fourier integrals in classical analysis}, volume 210 of {\em
  Cambridge Tracts in Mathematics}.
\newblock Cambridge University Press, Cambridge, second edition, 2017.

\bibitem{Wol03}
T.~H. Wolff.
\newblock {\em Lectures on harmonic analysis}, volume~29 of {\em University
  Lecture Series}.
\newblock American Mathematical Society, Providence, RI, 2003.
\newblock With a foreword by Charles Fefferman and preface by Izabella \L aba,
  Edited by \L aba and Carol Shubin.

\end{thebibliography}

\end{document}